\documentclass{conm-p-l}

\usepackage{amssymb}

\usepackage{graphicx}
\usepackage{color}

\usepackage{}

\newtheorem{theorem}{Theorem}[section]

\newtheorem{proposition}[theorem]{Proposition}

\theoremstyle{definition}
\newtheorem{definition}[theorem]{Definition}
\newtheorem{example}[theorem]{Example}

\theoremstyle{remark}
\newtheorem{remark}[theorem]{Remark}
\newtheorem{question}[theorem]{Question}

\numberwithin{equation}{section}

\newcommand{\K}{\Bbbk}
\newcommand{\MOn}{\overline{M}_{0,n}}

\newcommand{\PuiseuxC}{\ensuremath{\mathbb C \{\!\{t \} \! \}}}

\DeclareMathOperator{\relint}{relint}
\DeclareMathOperator{\trop}{trop}
\DeclareMathOperator{\inn}{in}
\DeclareMathOperator{\pos}{pos}
\DeclareMathOperator{\im}{im}

\DeclareMathOperator{\val}{val}
\DeclareMathOperator{\cl}{cl}
\DeclareMathOperator{\starr}{star}
\DeclareMathOperator{\spann}{span}
\DeclareMathOperator{\face}{face}
\DeclareMathOperator{\Aut}{Aut}

\begin{document}

\title{Introduction to tropical algebraic geometry}


\author{Diane Maclagan}
\address{Mathematics Institute, University of Warwick, Coventry, CV4 7AL, United Kingdom}
\email{D.Maclagan@warwick.ac.uk}
\thanks{Partially supported by  EPSRC grant EP/I008071/1}


\date{}

\begin{abstract}
This is an expository introduction to tropical algebraic geometry
based on my lectures at the Workshop on Tropical Geometry and
Integrable Systems in Glasgow, July 4--8, 2011, and at the ELGA 2011 school
on Algebraic Geometry and Applications in Buenos Aires, August 1--5,
2011.
\end{abstract}

\maketitle

\section{Introduction}

Tropical algebraic geometry is algebraic geometry over the tropical
semiring (Definition~\ref{d:tropicalsemiring}).  This replaces an
algebraic variety by a piecewise linear object which can be studied
using polyhedral combinatorics.

Tropical geometry has exploded as an area of research in the last
decade, with many new connections and applications appearing each
year.  These include enumerative geometry, mirror symmetry, arithmetic
geometry, and integrable systems. It builds on the older area of
tropical mathematics, more commonly known as max-plus algebra, which
arises in semigroup theory, computer science, and optimization.  The
name ``tropical'' was coined by some French mathematicians in honor of
the Brazilian computer scientist Imre Simon.  See \cite{ButkovicBook}
or \cite{GaubertSurvey} for an introduction to this older area.

The goal of this expository and elementary article is to introduce
this exciting new area.  We develop the theory of tropical varieties
and outline their structure and connection with ``classical''
varieties.

There are several approaches to tropical geometry.  We follow the
``embedded'' approach, which focuses on tropicalizing classical
varieties.  Another important branch of the subject focuses on
developing an abstract theory of tropical varieties in their own
right.  See work of Mikhalkin and collaborators \cite{MikhalkinIMPA},
\cite{MikhalkinICM} for details on this.  This direction is most
developed for curves.  One of the attractions of tropical geometry is
that it has so many disparate, but connected, facets, so any survey is
necessarily incomplete.

We begin by introducing tropical mathematics, focusing on tropical
polynomials and their solutions.  In this paper we will follow the
minimum convention for the tropical semiring:

\begin{definition}  \label{d:tropicalsemiring} The {\em tropical semiring} is $\mathbb R \cup \{ \infty \}$, 
with operation $\oplus$ and $\otimes$ given by $a \oplus b =
\min(a,b)$ and $a \otimes b = a+b$.  
\end{definition}

The tropical semiring is associative and distributive, with additive
identity $\infty$ and multiplicative identity $0$.  This satisfies
every axiom of a ring except additive inverses, so is a semiring.

Tropical operations are often simpler than regular operations.  For
example we have the ``Freshman's dream'': $(x \oplus y)^n = x^n \oplus  y^n$.  The
following examples illustrate that tropical polynomials are piecewise
linear functions.

\begin{example} \label{e:tropicalpolynomials}
\begin{enumerate}
\item The tropical polynomial $F(x)=-2\otimes x^3 \oplus -1 \otimes x^2 \oplus 1
  \otimes x \oplus 5$ is $\min(3x-2, 2x-1, x+1, 5)$ in regular
  arithmetic.  This is the piecewise linear function whose graph is
  shown in Figure~\ref{f:tropicalpolynomial}.

\begin{figure}
\center{\resizebox{6cm}{!}{\input{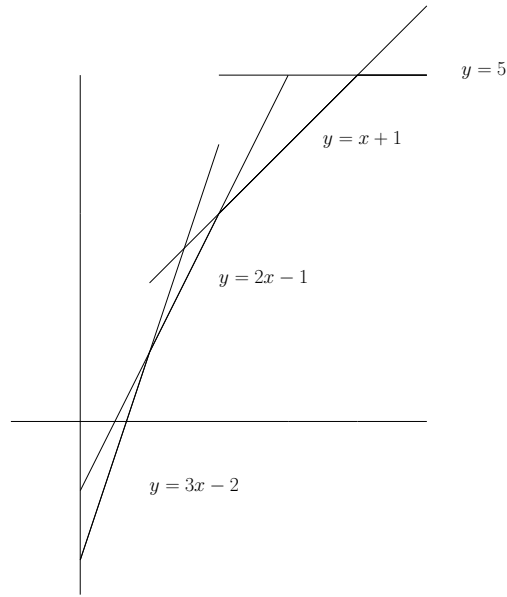}}}
\caption{A tropical polynomial \label{f:tropicalpolynomial}}
\end{figure}

\item \label{i:tropicalline} The tropical multivariate polynomial $x \oplus y \oplus 0$ is
  the piecewise linear function $\min(x,y,0)$ in regular arithmetic.
  Note that the zero cannot be removed here, as zero is not the
  additive identity.  This is a function from $\mathbb R^2$ to
  $\mathbb R$ whose domain is shown in Figure~\ref{f:tropicalline}.

\begin{figure}
\center{\resizebox{9cm}{!}{\input{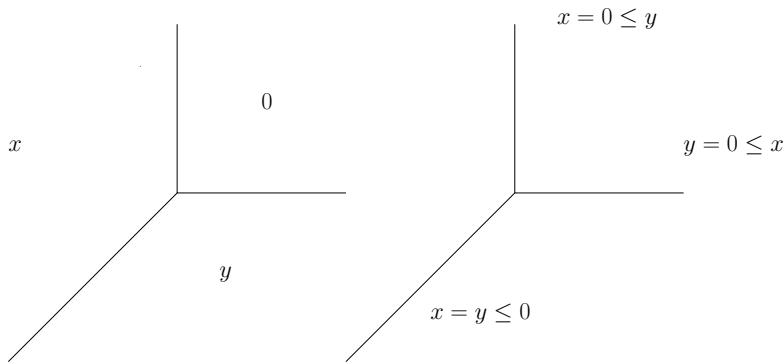}}}
\caption{A tropical line \label{f:tropicalline}}
\end{figure}

\end{enumerate}

\end{example}

With no subtraction, it is not obvious how to solve polynomial
equations.  For example, the equation $x \oplus 2 = 5$ has no
solution.  This problem has the following resolution.

\begin{definition}
The {\em tropical hypersurface} $V(F)$ defined by the tropical
polynomial $F$ in $n$ variables is the nondifferentiable locus of $F$
in $\mathbb R^n$.  This is the set of $x \in \mathbb R^n$ for which
the minimum is achieved at least twice in $F(x)$.
\end{definition}

\begin{example}
For the first polynomial of Example~\ref{e:tropicalpolynomials}, $V(F)
= \{1, 2,  4 \}$.  For the second, $V(F)$ is the union of the three rays
shown on the right in Figure~\ref{f:tropicalline}.
\end{example}

\begin{example}
The tropical quadratic formula is particularly simple.  If
$F(x)=a\otimes x^2 \oplus b \otimes x \oplus c$, then the graph of $F$
is shown in Figure~\ref{f:quadraticpoly}.  Note that there are two
cases, depending on the sign of the tropical discriminant $a+c-2b$.
If $2b\leq a+c$ then $V(F) = \{ c-b,b-a\}$.  If $2b \geq a+c$ then
$V(F) = \{ (c-a)/2 \}$.  Compare this with the usual quadratic
formula, and the usual discriminant.  Higher degree polynomials have 
similarly easy solutions.

\begin{figure}
\center{\resizebox{6cm}{!}{\input{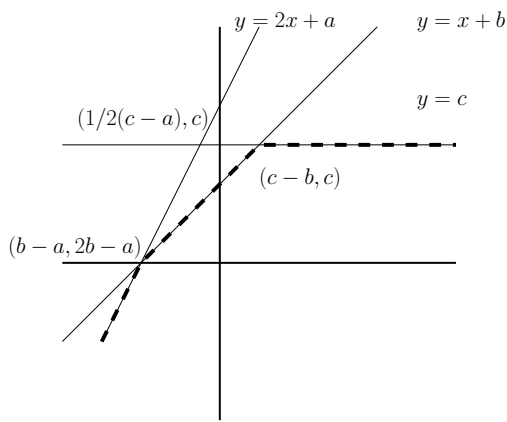}}}
\caption{Tropical quadratic polynomials \label{f:quadraticpoly}}
\end{figure}
\end{example}

\section{Tropical varieties}

The impact of tropical techniques comes when they are applied to
``classical'' objects, such as algebraic varieties.  We now explain
how to tropicalize certain algebraic varieties.

For a field $K$ we set $K^* = K \setminus \{0\}$.  Fix a valuation
$\val : K^* \rightarrow \mathbb R$.  This is a function 
satisfying:

\begin{enumerate}
\item $\val(ab)=\val(a)+\val(b)$
\item $\val(a+b) \geq \min(\val(a),\val(b))$.
\end{enumerate}

\begin{example}

\begin{enumerate}

\item $K=\mathbb C$ with the trivial valuation $\val(a)=0$ for all $a
  \in \mathbb C^*$.
\item $K = \PuiseuxC$, the field of Puiseux series.  This is the union $\bigcup_{n
  \geq 1} \mathbb C (( t^{1/n}))$, and is the algebraic closure of the
  field of Laurent series.  Elements are Laurent series with rational
  exponents where in any given series the exponents all have a common
  denominator.  The valuation of $a \in \PuiseuxC$ is
  the lowest exponent appearing.  For example,
  $\val(3t^{-1/2}+8t^2+7t^{13/3}+ \dots)=-1/2$.

\item $K=\mathbb Q$ or $\mathbb Q_p$ with the $p$-adic valuation.  If
  $q=p^na/b \in \mathbb Q$ with $p$ not dividing $a$ or $b$, then $\val_p(q) =n$.
  For example, $\val_2(8)=3$, and $\val_3(5/6)=-1$.

\end{enumerate}
\end{example}

\begin{definition}
The {\em tropicalization} of a Laurent polynomial $f = \sum c_u x^u \in
K[x_1^{\pm 1},\dots,x_n^{\pm 1}]$ is $\trop(f) : \mathbb R^n
  \rightarrow \mathbb R$ given by
$$\trop(f)(w) = \min(\val(c_u) + w \cdot u).$$ This is obtained by
  regarding the addition and multiplication as tropical addition and
  multiplication, and changing the coefficients to their valuations.
\end{definition}

\begin{example}
Let $K=\mathbb Q$ with the $2$-adic valuation, and let $f =
6x^2+5xy+10y^2+3x-y+4 \in \mathbb Q[x^{\pm 1}, y^{\pm 1}]$.  Then
$\trop(f) = \min(2x+1,x+y,2y+1,x,y,2)$.  This is illustrated in
Figure~\ref{f:tropicalquadric}.

\begin{figure}
\center{\resizebox{8cm}{!}{\input{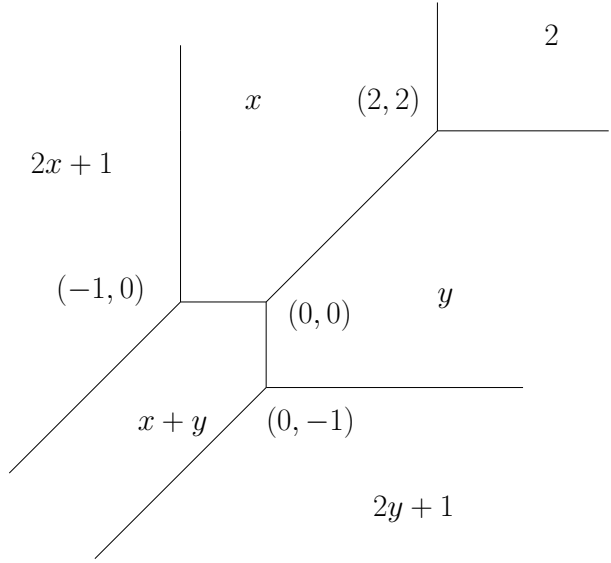}}}
\caption{A tropical quadric \label{f:tropicalquadric}}

\end{figure}

\end{example}

For $f \in K[x_1^{\pm 1}, \dots, x_n^{\pm 1}]$ the (classical)
hypersurface $V(f)$ equals $\{x \in (K^*)^n : f(x)=0 \}$.  The 
  tropicalization $\trop(V(f))$ of $V(f)$ is the tropical
hypersurface of $\trop(f)$.  This is the nondifferentiable locus of
$\trop(f)$, or equivalently:
$$\trop(V(f)) = \{ w \in \mathbb R^n : \text{ the minimum in } \trop(f)(w) \text{ is achieved at least twice} \}.$$

Note that $\trop(x^uf) = \trop(x^u) + \trop(f)$, so $\trop(x^uf)(w) =
w \cdot u +\trop(f)(w)$, and thus $\trop(V(x^uf))=\trop(V(f))$.  This
explains why the natural algebraic varieties to tropicalize are
subvarieties of $T=(K^*)^n$, rather than subvarieties of $\mathbb A^n$
or $\mathbb P^n$, as monomial functions are invertible on $T$.

Let $Y=V(I)$ be a subvariety of $T=(K^*)^n$. If $I= \langle f_1,\dots,f_r \rangle$ then 
\begin{align*}
Y =V(I) & = \{ x \in T : f_1(x) = \dots = f_r(x) =0 \}\\
& = \{ x \in T : f(x) = 0 \text{ for all } f \in I \}\\
& = \bigcap_{f \in I} V(f).\\
\end{align*}

\begin{definition} \label{d:tropY}
The tropicalization of a variety $Y=V(I)  \subseteq T$ is 
$$\trop(Y) = \bigcap_{f \in I} \trop(V(f)).$$
\end{definition}

\begin{example}
\label{eg:tropicalvariety}

\begin{enumerate}

\item  \label{e:lineality} $Y = V(x+y+z+w,x+2y+5z+11w) \subseteq (\mathbb C^*)^4$.  Then
  $\trop(Y)$ has the property that if $\alpha \in \trop(Y)$, then
  $\alpha +\lambda(1,1,1,1) \in \trop(Y)$ for all $\lambda \in \mathbb R$.
  We can thus quotient by the span of $(1,1,1,1)$ and describe
  $\trop(Y)$ in $\mathbb R^4/\mathbb R(1,1,1,1) \cong \mathbb R^3$.
  This consists of four rays, being the images of the positive
  coordinate directions $(1,0,0,0)$, $(0,1,0,0)$, $(0,0,1,0)$, and
  $(0,0,0,1)$.

\item Let $Y= V(t^3x^3+x^2y+xy^2+t^3y^3+x^2+t^{-1}xy+y^2+x+y+t^3)
  \subseteq (\PuiseuxC^*)^2$.  Then $\trop(V(f))$ is shown in
  Figure~\ref{f:degree3}.  This is a ``tropical elliptic curve''; see,
  for example, \cite{Vigeland} or \cite{KatzMarkwigMarkwig}.

\begin{figure}
\center{\resizebox{7cm}{!}{\input{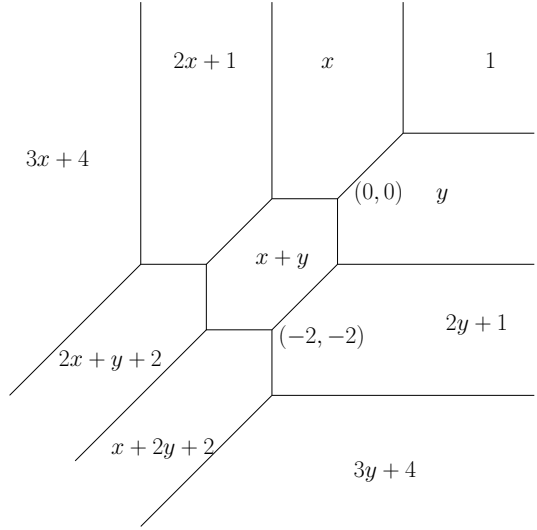}}}
\caption{A tropical elliptic curve \label{f:degree3}}
\end{figure}

\end{enumerate}

\end{example}

\begin{remark}  If $Y=V(f_1,\dots,f_r)$, the tropical variety $\trop(Y)$ does not always
equal $\bigcap_{i=1}^r \trop(V(f_i))$.  For example, in the first part of
Example~\ref{eg:tropicalvariety},  $\trop(V(x+y+z+w))=\trop(V(x+2y+5z+11w))$, but this 
is larger than the tropical variety $\trop(Y)$.
\end{remark}

\begin{definition}
If $Y =V(I)$, a set $\{f_1,\dots,f_r \} \subseteq  I$ with 
$$\trop(Y) = \bigcap_{i=1}^r \trop(V(f_i))$$ is called a {\em tropical
basis} for $I$.  Finite tropical bases always exist (see
\cite[Theorem 11]{BJSSTCompute}), so $\trop(Y)$ is a piecewise linear object.
\end{definition}

The following is the guiding question behind much tropical geometry research:

\begin{question} \label{q:invariants}
Which properties of $Y$ or of
compactifications of $Y$ can be deduced from $\trop(Y)$?
\end{question}

\section{Drawing curves in the plane}

Before delving deeper into the connections between classical and
tropical varieties, with a view to partial answers to
Question~\ref{q:invariants}, we first describe how to draw plane
curves.

Let $C=V(f) \subseteq (K^*)^2$ for $f \in K[x^{\pm 1},y^{\pm 1}]$.  If
$f = \sum c_{ij} x^iy^j$, then $\trop(f) = \min(\val(c_{ij})+ ix+jy)$,
and $\trop(C) = \trop(V(f))$ is the locus where the minimum is achieved
at least twice.  Some elementary polyhedral geometry gives a short-cut
to compute this locus, so we first review this notation.  An excellent general introduction to polyhedral geometry is~\cite{Ziegler}.

\begin{definition}
A {\em polyhedron} in $\mathbb R^n$ is the intersection of finitely many
half-spaces in $\mathbb R^n$.  This can be written as:
$$P = \{ x \in \mathbb R^n : Ax \leq b \},$$ where $A$ is a $d \times
n$ matrix, and $b \in \mathbb R^d$.  The dimension of $P$ is the
dimension of the subspace $\ker(A)$.  For a subgroup $\Gamma \subseteq
\mathbb R$, we say that a polyhedron $P$ is $\Gamma$-rational if $A$
has rational entries, and $b \in \Gamma^d$.  When $\Gamma = \mathbb
Q$, we say that the  polyhedron is rational.  

If $b=0$, then $P$ is called a {\em cone}.  In that case there are
$\mathbf{v}_1,\dots,\mathbf{v}_s$ for which $P =
\pos(\mathbf{v}_1,\dots,\mathbf{v}_s) := \{ \sum_{i=1}^s \lambda_i
\mathbf{v}_i : \lambda_i \geq 0 \}$.

The {\em face} of a polyhedron
$P$ determined by $w \in (\mathbb R^n)^*$ is the set
$$\face_w(P) = \{ x \in P : w \cdot x \leq w \cdot y \text{ for all } y \in P\}.$$
\end{definition}

\begin{example}
Let $P \subset \mathbb R^2$ be the square with vertices $\{(0,0), (1,0), (0,1), (1,1) \}$.  This has the description 
$$P = \left\{ x \in \mathbb R^2 : \left( \begin{array}{rr} 1 & 0 \\ 0 & 1
  \\ -1 & 0 \\ 0 & -1 \\ \end{array} \right) x \leq 
\left( \begin{array}{r} 1 \\ 1\\ 0 \\ 0\\ \end{array} \right) \right\}.$$

Then we have
\begin{enumerate}
\item $\face_{(1,0)}(P)$ is the edge of the square with vertices $\{(0,0), (0,1) \}$,
\item $\face_{(1,1)}(P)$ is the vertex $(0,0)$, and
\item $\face_{(0,0)}(P)$ is $P$.
\end{enumerate}

\end{example}

\begin{definition}
A {\em polyhedral complex} $\Sigma$ is a finite union of polyhedra for which
any nonempty intersection of two polyhedra $\sigma_1, \sigma_2 \in
\Sigma$ is a face of each.  If every polyhedron in $\Sigma$ is a cone,
then $\Sigma$ is called a fan.  

The {\em normal fan} to a polyhedron $P$ is the fan $N(P)$ with cones
$C[w] = \cl(w' : \face_{w'}(P) = \face_w(P) )$, where $\cl( \cdot )$ is
the closure in the usual Euclidean topology on $\mathbb R^n$.  This is
sometimes called the inner normal fan, as the cones corresponding to
facets (faces of maximal dimension) are the inner normal vectors to
these faces
\end{definition}

\begin{definition}
Fix $f \in K[x_1^{\pm 1},\dots, x_n^{\pm 1}]$.  The {\em Newton polytope} of
$f = \sum c_u x^u $ is the convex hull of the $u \in \mathbb Z^n$ with $c_u \neq 0$:
$P = \{ \sum_{c_u \neq 0}  \lambda_u u : \sum \lambda_u =1 \}$.  
\end{definition}

\begin{example}
Let $f = x^2y+5y^2-3x+2$.  The Newton polytope $P$ of $f$ and its
normal fan $N(P)$ are illustrated in Figure~\ref{f:normalfan}.  The
polytope $P$ is the convex hull of the exponent vectors $\{ (2,1),
(0,2), (1,0), (0,0) \}$.

\begin{figure}[h]
\center{\resizebox{!}{5cm}{\includegraphics{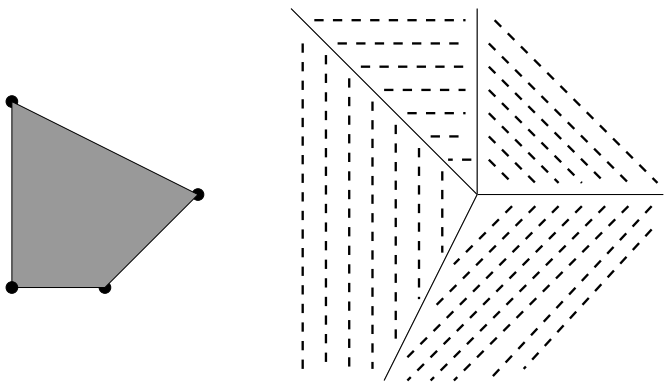}}}
\caption{Newton polytope and normal fan \label{f:normalfan}}

\end{figure}

\end{example}

We now temporarily restrict our attention to plane curves, so $n=2$.
The tropicalization of a plane curve $C=V(f)$ for $f = \sum c_{ij}x^iy^j$ depends
on the valuations of the coefficients $c_{ij}$.  If the valuation on $K$
is trivial ($\val(a)=0$ for all $a \neq 0$), then $\trop(C)$ is the
union of all one-dimensional cones in the normal fan $N(P)$ to the
Newton polytope of $f$.

We now consider the case that $K$ has a nontrivial valuation.  Let $\widetilde{P}$ be
the convex hull of the set $\{ (i,j, \val(c_{ij})) : c_{ij} \neq 0 \}$ in
$\mathbb R^{2+1}=\mathbb R^3$, and let $N(\widetilde{P})$ be its
normal fan.  The regular subdivision $\Delta_{(\val(c_{ij}))}$ of $P$
corresponding to the vector $(\val(c_{ij}))$ is the projection to $P$ of
the ``lower faces'' of $\widetilde{P}$.  More information about regular subdivisions can be found in \cite[Chapter 7]{GKZ} and \cite[Chapter 2]{TriangulationsBook}.

\begin{example}
Let $f = 2x^2+xy-6y^2+5x-3y+2 \in \mathbb Q[x^{\pm 1}, y^{\pm 1}]$
where $\mathbb Q$ has the $2$-adic valuation.  Then the regular
subdivision of the Newton polytope of $f$ corresponding to $\val(c_u)$
is shown in Figure~\ref{f:regsubdivision}.

\begin{figure}[h]

\center{\resizebox{!}{3cm}{\includegraphics{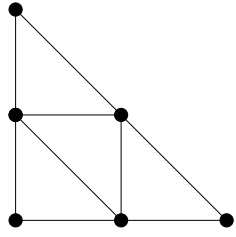}}}
\caption{\label{f:regsubdivision} The regular subdivision induced by $(\val(c_u))$}

\end{figure}

\end{example}

The tropicalization of $C$ is  $\trop(C) = \{ w \in \mathbb R^2 :
\face_{(w,1)}(\widetilde{P}) \text{ is not a vertex } \}$.  This is
the reflection of the dual graph to $\Delta_{(\val(c_{ij}))}$ under $x
\mapsto -x$.

\begin{example}

Let $f =27x^3+6x^2y+12xy^2+81y^2+3x^2+5xy+3y^2+3x+2y+243 \in
  \mathbb Q[x^{\pm 1},y^{\pm 1}]$ where $\mathbb Q$ has the $3$-adic
  valuation.

Then the regular triangulation corresponding to $(\val(c_{ij}))$ is shown in Figure~\ref{f:cubic}, along with the tropical variety.

\begin{figure}[h]
\center{\resizebox{!}{3cm}{\includegraphics{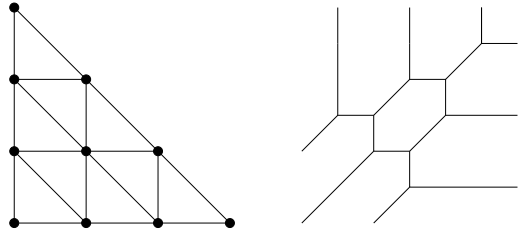}}}
\caption{A tropical plane cubic \label{f:cubic}}

\end{figure}

\end{example}

\section{The Fundamental and Structure theorems}

The definition of the tropical variety in Definition~\ref{d:tropY} was
by analogy with the classical case.  For an ideal $I$ the tropical
variety $\trop(V(I))$ is the set of common ``tropical zeros'' of the
tropicalizations of the polynomials $f \in I$.  The connection between
the tropicalization $\trop(Y)$ of a variety $Y=V(I) \subseteq (K^*)^n$
and the original variety $Y$ is closer than this analogy might
suggest, as the following theorem shows.

\begin{theorem}[Fundamental theorem of tropical algebraic geometry]
\label{t:fundamentaltheorem}
Let $K$ be an algebraically closed field with a nontrivial valuation
$\val: K^* \rightarrow \mathbb R$, and let $Y$ be a subvariety of
$(K^*)^n$.  Then
\begin{align*} \trop(Y) & = \cl(\val(Y))\\
& = \cl( (\val(y_1),\dots,\val(y_n)) : y=(y_1,\dots,y_n) \in Y),\\
\end{align*}
where the closure is in the usual Euclidean topology on $\mathbb R^n$.
\end{theorem}

\begin{example}
Let $Y = V(x+y+1) \subseteq (K^*)^2$, where $K=\PuiseuxC$.  Then $Y = \{ (a,-1-a) : a \in K^* \setminus \{-1\} \}$.  Note that 

$$( \val(a), \val(-1-a) ) =  
\begin{cases}
(\val(a), 0)  & \text{if }\val(a) >0, \\
(\val(a), \val(a) )  & \text{if }\val(a)<0, \\
(0,\val(b))  & \text{if }a = -1 +b, \val(b)>0, \\
(0,0)  & \text{ otherwise}.\\
\end{cases}.
$$

Note that, as predicted by Theorem~\ref{t:fundamentaltheorem}, the
union of these sets is precisely $\trop(Y)$, as calculated in
part~\ref{i:tropicalline} of Example~\ref{e:tropicalpolynomials}.
This is illustrated in Figure~\ref{f:tropicallineValuations}, which
should be compared with Figure~\ref{f:tropicalline}.

\begin{figure}
\center{\resizebox{4cm}{!}{\input{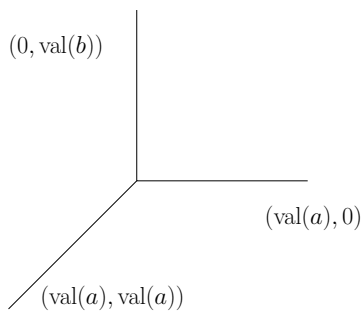}}}
\caption{The Fundamental Theorem applied to a tropical line \label{f:tropicallineValuations}}
\end{figure}

\end{example}

\begin{example}
Let $f= 4x^2+xy-4y^2+x-y-4 \in \mathbb Q[x,y]$, where $\mathbb Q$ has
the $2$-adic valuation.  Then $\trop(V(f))$ is shown in
Figure~\ref{f:FT2adic}.  Note that the point $(2,2) \in V(f)$, and
$(\val(2),\val(2)) = (1,1) \in \trop(V(f))$.

\begin{figure}
\center{\resizebox{!}{6cm}{\input{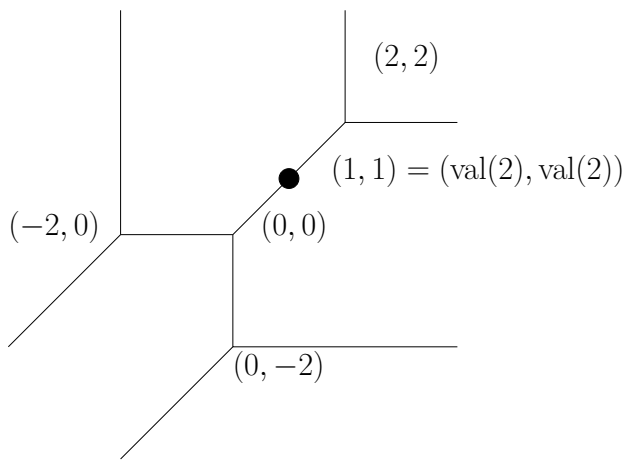}}}
\caption{\label{f:FT2adic} An example of the Fundamental Theorem}
\end{figure}

\end{example}

Theorem~\ref{t:fundamentaltheorem} was first shown for hypersurfaces
by Kapranov; see~\cite{EKL}.  The general case first appears in work of 
Speyer and Sturmfels~\cite{SpeyerThesis},
\cite{SpeyerSturmfelsTropicalGrassmannian}; see also the work of
Draisma~\cite{Draisma}, Payne~\cite{PayneFiber}, and Jensen, Markwig,
and Markwig~\cite{JensenMarkwigMarkwig}.  The hard part of
Theorem~\ref{t:fundamentaltheorem} is to show that if $w \in \trop(Y)
\cap (\im \val)^n$ then there is $y \in Y$ with $\val(y)=w$.  Showing
that $\{ \val(y) : y \in Y \}  \subseteq \trop(Y)$ is comparatively
easy.

The slogan form of Theorem~\ref{t:fundamentaltheorem} is then:

\begin{quotation}
{\em Tropical varieties are combinatorial shadows of classical varieties.}
\end{quotation}

The word ``combinatorial'' is justified by the Structure Theorem for
tropical varieties, which gives combinatorial constraints on which
sets can be tropical varieties.  The statement of this theorem
requires some more polyhedral definitions.

\begin{definition}
Let $\Sigma$ be a polyhedral complex.  The {\em support} of $\Sigma$ is the
set
$$|\Sigma| = \{ x \in \mathbb R^n : x \in \sigma \text{ for some } \sigma \in \Sigma \}.$$

A polyhedral complex is {\em pure} if the dimension of every maximal
polyhedron is the same.  See Figure~\ref{f:pure}.

\begin{figure}
\center{\resizebox{!}{3cm}{\includegraphics{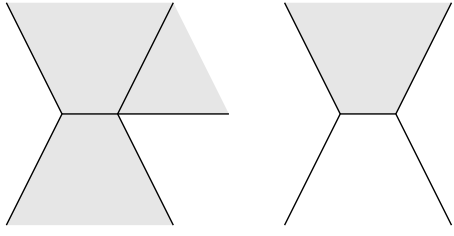}}}
\caption{\label{f:pure} The complex on the left is pure, while the one
  on the right is not.}
\end{figure}

The {\em lineality space} $L$ of a polyhedral complex $\Sigma$ is the
largest subspace of $\mathbb R^n$ for which $x +l \in \Sigma$ for all
$x \in \Sigma, l \in L$.  In part~\ref{e:lineality} of Example~\ref{eg:tropicalvariety} the lineality
space of $\trop(Y)$ is $\spann((1,1,1,1)) \subseteq \mathbb R^4$.
\end{definition}

\begin{definition}
A {\em weighted polyhedral complex} is a pure polyhedral complex $\Sigma$ with a
weight $w_{\sigma} \in \mathbb N$ for all maximal-dimensional $\sigma
\in \Sigma$.

Let $\Sigma$ be a weighted $(\im \val)$-rational polyhedral complex that
is pure of dimension $d$.  The complex $\Sigma$ is {\em balanced} if
the following ``zero-tension'' conditions hold.

\begin{enumerate}
\item If $\Sigma$ is a one-dimensional rational fan, let
  $\mathbf{u}_1,\dots,\mathbf{u}_s$ be the first lattice points on the
  rays of $\Sigma$, and let $w_i$ be the weight of the cone containing
  the lattice point $\mathbf{u}_i$.  Then $\Sigma$ is balanced if
  $\sum_{i=1}^s w_i \mathbf{u}_i =0$. 

\item For a general polyhedral complex $\Sigma$, fix a $(d-1)$-dimensional polyhedron $\tau$
  of $\Sigma$.  Let $L = \spann( x-y : x, y \in \tau)$ be the affine
  span of $\tau$.  Let $\starr_{\Sigma}(\tau)$ be the rational
  polyhedral fan whose support is $\{w \in \mathbb R^n : \text{ there
    exists } \epsilon >0 \text{ for which } w' + \epsilon w \in \Sigma
  \text{ for all } w' \in \tau \}+ L$.  This has one cone for each
  polyhedron $\sigma \in \Sigma$ that contains $\tau$, and has
  lineality space $L$.  The quotient $\starr_{\Sigma}(\tau)/L$ is a
  one dimensional fan which inherits weights from $\Sigma$.  We say
  that $\Sigma$ is balanced at $\tau$ if the one-dimensional fan
  $\starr_{\Sigma}(\tau)/L$ is balanced.  The polyhedral complex
  $\Sigma$ is balanced if $\Sigma$ is balanced at all
  $(d-1)$-dimensional cones.
\end{enumerate}

\end{definition}

\begin{example}
Let $f=x^2y^2+x^3+y^3+1 \in \mathbb C[x^{\pm 1},y^{\pm 1}]$.  Then
$\trop(V(f))$ is a one-dimensional fan with four rays: $\pos((1,0))$,
$\pos((0,1))$, $\pos((-2,-1))$, and $\pos((-1,-2))$.  This is balanced with weights
$3,3,1$, and $1$:
$3(1,0)+3(0,1)+1(-2,-1)+1(-1-2)=(0,0)$.
\end{example}

\begin{definition}
We associate a graph to a pure polyhedral complex $\Sigma$, with 
a vertex for each maximal polyhedron $\sigma \in \Sigma$, and
an edge between two vertices if the corresponding polyhedra intersect
in a codimension-one face.  The polyhedral complex $\Sigma$ is {\em connected
  through codimension-one} if this graph is connected.  For example, the
polyhedral complex on the left of Figure~\ref{f:connected} is
connected through codimension-one, while the one on the right is not.

\begin{figure}
\center{\resizebox{!}{4cm}{\input{connectedincodimone.pstex_t}}}
\caption{A polyhedral complex that is (left) and is not (right) connected through codimension one \label{f:connected}}
\end{figure}

\end{definition}

Recall that a variety $Y \subseteq (K^*)^n$ is irreducible if we
cannot write $Y = Y_1 \cup Y_2$ for $Y_1, Y_2$ nonempty proper 
subvarieties of $Y$.  Note that by Theorem~\ref{t:fundamentaltheorem}
we have $\trop(Y_1 \cup Y_2) = \trop(Y_1) \cup \trop(Y_2)$.

The following ``Structure Theorem'' summarizes the combinatorial
structure on the tropical variety.  This is essentially due to Bieri
and Groves~\cite{BieriGroves}, with the stronger connectedness
statement from \cite{BJSSTCompute} and \cite{CartwrightPayne}.

\begin{theorem}[Structure Theorem] \label{t:structure}
Let $Y$ be a $d$-dimensional irreducible subvariety of $(K^*)^n$.
Then $\trop(Y)$ is the support of a pure $d$-dimensional weighted
balanced $(\im \val)$-rational polyhedral complex that is connected
through codimension-one.
\end{theorem}

The Structure Theorem means that the tropical variety has a discrete
structure, and records information about the original variety (such as
its dimension).

The weights $w_{\sigma}$ on maximal cones that make $\trop(Y)$
balanced can be computed from $Y$, as we explain in the next section.

\section{The computational approach to tropical varieties}

An important aspect of tropical varieties is that they can actually be
computed in practice.  This uses an extension of the theory of
Gr\"obner bases to fields with a valuation.

Fix a splitting $(\im \val) \rightarrow K^*$ of the valuation.  This is
a group homomorphism $u \mapsto t^u$ with $\val(t^u)=u$.  For example,
when $K= \mathbb C$ with the trivial valuation ($\val(a)=0$ for all $a
\neq 0$), then we can choose the splitting $0 \mapsto 1$.  When
$K=\PuiseuxC$, we can choose $u \mapsto t^u$, and when $K=\mathbb Q$
with the $p$-adic valuation, we can choose $u \mapsto p^u$.  Such
splittings always exist when $K$ is algebraically closed; see
\cite[Lemma 2.1.13]{TropicalBook}.

Let $R = \{a \in K : \val(a) \geq 0 \}$ be the valuation ring of $K$.
The ring $R$ is local, with maximal ideal $\mathfrak{m} = \{a \in K :
\val(a) >0 \} \cup \{0 \}$.  The quotient $\K = R/\mathfrak{m}$ is the
residue field.

\begin{example}
\begin{enumerate}
\item 
When $K=\mathbb C$ has the trivial valuation, we have $R=\mathbb C$,
and $\mathfrak{m}=0$, so $\K=\mathbb C$.
\item When $K = \PuiseuxC$, $R=\bigcup_{n \geq 1} \mathbb C[\![t^{1/n}]\!]$, and
  $\K=\mathbb C$.
\item When $K=\mathbb Q_p$, $R=\mathbb Z_p$, and $\K = \mathbb
  Z/p\mathbb Z$.
\end{enumerate}

\end{example}

Given a polynomial $f = \sum c_u x^u \in K[x_1^{\pm 1},\dots,x_n^{\pm 1}]$, 
and $w \in (\im \val)^n$, the {\em initial form} is 
$$\inn_w(f) = \sum_{\val(c_u)+w \cdot u = \trop(f)(w)}
\overline{t^{-\val(c_u)}c_u} x^u \in \K[x_1^{\pm 1},\dots,x_n^{\pm
    1}],$$ where for $a \in R$ we denote by $\overline{a}$ the image
of $a$ in $\K$.

\begin{example}

Let $K=\mathbb Q$ with the $2$-adic valuation, and let
$f=2x^2+xy+6y^2+5x-3y+4$.  Then for $w=(2,2)$ we have $\trop(f)(w)=2$,
so $\inn_w(f) = \overline{5}x+\overline{-3}y+\overline{2^{-2}4} = x+y+1
\in \mathbb Z/2\mathbb Z[x^{\pm 1},y^{\pm 1}]$.

For $w =(-2,-1)$ we have $\trop(f)(w) = -3$, so $\inn_w(f) = x^2+xy$.

\end{example}

Given an ideal $I \subset K[x_1^{\pm 1},\dots,x_n^{\pm 1}]$ and $w \in (\im \val)^n$, the initial ideal of $I$ is 
$$\inn_w(I) = \langle \inn_w(f) : f\in I \rangle.$$

This is a variant of the standard theory of Gr\"obner bases for ideals
in a polynomial ring; see~\cite{CLO} for an excellent introduction to the 
standard case.  As in that setting, the initial ideal need not be
generated by the initial forms of a generating set for $I$.  There are
always finite generating sets for $I$, however, for which this is the
case.  These finite sets (Gr\"obner bases) can be computed using a
variant of the standard Gr\"obner basis algorithm.

The connection of this theory of Gr\"obner bases to tropical geometry
is the following computational characterization of tropical varieties,
due to Sturmfels~\cite[Chapter 9]{SolvingPolynomialEquations}.

\begin{proposition} \label{p:tropicalGrobner}
Let $Y = V(I) \subseteq (K^*)^n$ and let $w \in (\im \val)^n$.  Then
$w \in \trop(Y)$ if and only if $\inn_w(I) \neq  \langle 1 \rangle
\subseteq \K[x_1^{\pm 1},\dots,x_n^{\pm 1}]$.
\end{proposition}

\begin{example} \label{e:Grobnercomplex}
Let $I = \langle 2x^2+xy+6y^2+5x-3y+4 \rangle \in \mathbb Q[x^{\pm
    1},y^{\pm 1}]$, where $\mathbb Q$ has the $2$-adic valuation.  The
claim of Proposition~\ref{p:tropicalGrobner} is illustrated in
Figure~\ref{f:GrobnerComplex}. 

\begin{figure}
\center{\resizebox{!}{6cm}{\input{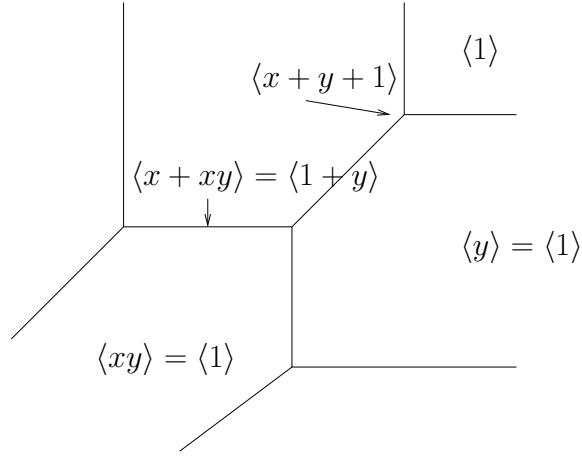}}}
\caption{\label{f:GrobnerComplex}The Gr\"obner complex of Example~\ref{e:Grobnercomplex}}

\end{figure}

\end{example}

There is a polyhedral complex $\Sigma$ with $\inn_w(I)$ constant for
$w \in \relint(\sigma)$ for any $\sigma \in \Sigma$.  This is the {\em
  Gr\"obner complex} of a homogenization of $I$.  In the case that
$K=\mathbb C$ with the trivial valuation the Gr\"obner complex is a
polyhedral fan, called the Gr\"obner fan.  See \cite{BayerMorrison},
\cite{MoraRobbiano}, \cite[Chapter 2]{GBCP} for more on the Gr\"obner
fan, and \cite{Bellairs}, \cite{TropicalBook} for expositions of the
Gr\"obner complex.

Proposition~\ref{p:tropicalGrobner}
implies that $\trop(Y)$ is the union of the polyhedra $\sigma$ in the
Gr\"obner complex of $I(Y)$ for which $\inn_w(I) \neq \langle 1
\rangle$ for any $w \in \relint(\sigma)$.
The software {\tt{gfan}} \cite{gfan} by Anders Jensen computes tropical
varieties by exploiting this Gr\"obner description.

The Structure Theorem (Theorem~\ref{t:structure}) asserts that
$\trop(Y)$ is the support of a weighted balanced polyhedral complex.
The weights that make the polyhedral complex balanced can be defined
using the Gr\"obner theory we have just summarized.  We assume that
the polyhedral complex structure $\Sigma$ has been chosen, as
described above, so that $\inn_w(I(Y))$ is constant on the relative
interior of polyhedra of $\Sigma$.
  Fix $w \in (\im \val)^n$ in the
relative interior of a maximal polyhedron of $\Sigma$.  Then
$V(\inn_w(I(Y))) \subseteq (\K^*)^n$ is a union of $(\K^*)^d$-orbits,
where $d=\dim(Y)$ (see \cite[Chapter 3]{TropicalBook} for details).
We set $w_{\sigma}$ to be the number of such orbits (counted with
multiplicity).

Hidden in the proof of the Structure Theorem is the fact that this
choice makes the polyhedral complex $\Sigma$ balanced.

\begin{example}
Let $f = x^2+3x+2+x^2y+2xy^2-2y^2 \subseteq \mathbb C[x^{\pm 1},y^{\pm
    1}]$.  Then $\trop(V(f))$ is a one-dimensional fan with five rays,
spanned by the vectors $\{(1,0), (0,1), (-1,0), (-1,-1), (0,-1) \}$.
When $w=(0,1)$, $\trop(f)(w)=0$, so $\inn_w(f) = x^2+3x+2=(x+2)(x+1)$.
We then have $V(\inn_w(f)) = \{ (-2,a) : a \in \mathbb C^* \} \cup \{ (-1,a) :
a \in \mathbb C^* \}$. The weight on the cone spanned by $(0,1)$ is
thus $2$.  When $w = (-1,-1)$, $\inn_w(f) = x^2y+2xy^2$, so
$V(\inn_w(f)) = V(x+2y) = \{ (2a,-a) : a \in \mathbb C^* \}$.  This
makes the weight on this cone $1$.  Similarly, the weight on the cone
spanned by $(1,0)$ is $2$ and all other weights are $1$.  Note that
$$2\left( \begin{array}{r} 1 \\ 0 \end{array} \right) +
2\left( \begin{array}{r} 0 \\ 1 \end{array} \right) +
1\left( \begin{array}{r} -1 \\ 0 \end{array} \right) +
1\left( \begin{array}{r} -1 \\ -1 \end{array} \right) +
1\left( \begin{array}{r} 0 \\ -1 \end{array} \right) =
\left( \begin{array}{r} 0 \\ 0 \\ \end{array} \right).$$

This is illustrated in Figure~\ref{f:balancing}.

\begin{figure}
\center{\resizebox{!}{4cm}{\input{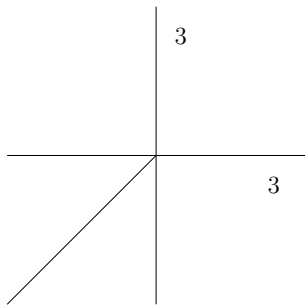}}}
\caption{The balancing condition \label{f:balancing}}

\end{figure}

\end{example}

\section{Connections to subvarieties of toric varieties}

Tropical varieties record information not just about subvarieties of
the algebraic torus $(K^*)^n$, but also about their compactifications.
This vastly extends the reach of tropical techniques; for example,
every projective variety is a compactification of a subvariety of the
algebraic torus.  The connection is through the theory of toric
varieties.

  We assume here that $K=\mathbb C$.  This means that for $Y \in
  (\mathbb C^*)^n$, $\trop(Y)$ is the support of a weighted balanced
  rational polyhedral fan.

\begin{definition}
A (normal) {\em toric variety} is a normal variety $X$ containing a dense
copy of $T=(\mathbb C^*)^n$ with an action of $T$ on $X$ that extends
the action of $T$ on itself.
\end{definition}

Examples include:
\begin{enumerate}
\item $X=(\mathbb C^*)^n$,

\item $X=\mathbb A^n \supset (\mathbb C^*)^n = \{x \in \mathbb A^n :
  x_i \neq 0 \text{ for } 1 \leq i \leq n \}$,

\item $X=\mathbb P^n \supset (\mathbb C^*)^n = \{x \in \mathbb P^n :
  x_i \neq 0 \text{ for } 0 \leq i \leq n \}$,

\item $X = \mathbb P^1 \times \mathbb P^1 \supset (\mathbb C^*)^2$.
\end{enumerate}

A toric variety $X$ is a union of $T$-orbits.  These can be recorded
using a polyhedral fan $\Sigma$.

\begin{example}
The projective plane $\mathbb P^2$ decomposes into the following
$T=(\mathbb C^*)^2$-orbits:
$$(\mathbb C^*)^2 \cup \{ [0: a:b] : a,b \in \mathbb C^* \} \cup \{
[a: 0:b] : a,b \in \mathbb C^* \} \cup \{ [a: b:0] : a,b \in \mathbb
C^* \} $$ $$\cup \{[1:0:0] \} \cup \{ [0:1:0] \} \cup \{0:0:1]\}.$$ The
corresponding fan is shown in Figure~\ref{f:P2orbits}.

\begin{figure}[h]
\center{\resizebox{!}{3cm}{\includegraphics{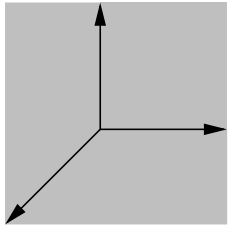}}}
\caption{\label{f:P2orbits} The fan of $\mathbb P^2$}
\end{figure}

\end{example}

Alternatively (and more standardly), given a rational polyhedral fan
$\Sigma$ we construct a toric variety $X_{\Sigma}$ by gluing together
torus orbits closures.  Each cone of $\Sigma$ determines an affine toric
variety, and the fan tells us how to glue them together.  For example,
for $\mathbb P^2$ the fan tells us to construct $\mathbb P^2$ by
gluing together the three affine charts $\{x \in \mathbb P^2 : x_i
\neq 0 \}$ for $0 \leq i \leq 2$.
 For more background on toric varieties, some good references include
 \cite{CLS} and \cite{fulton}.

The connection to tropical geometry begins with the following
question.

\begin{question}
Given a toric variety $X_{\Sigma}$, and a subvariety $\overline{Y} \subseteq X_{\Sigma}$, which $T$-orbits of $X_{\Sigma}$ does $\overline{Y}$ intersect?
\end{question}

Surprisingly, the answer uses tropical geometry.  The subvariety
$\overline{Y}$ intersects the torus orbit indexed by $\sigma \in
\Sigma$ if and only if $\trop(\overline{Y} \cap T)$ intersects
$\mathrm{relint}(\sigma)$.  This follows from work of Tevelev
\cite{TevelevCompactifications}.

\begin{example}
Let $\overline{Y}=V(x+y+z) \subseteq \mathbb P^2$.  Then $\trop(\overline{Y}
\cap T )= \trop(x+y+1)$, which is the standard tropical line.  This
intersects every one of the cones of the fan of $\mathbb P^2$ except for the
top-dimensional ones.  Indeed, the top-dimensional cones correspond to
the $T$-fixed points $[1:0:0], [0:1:0]$ and $[0:0:1]$, which do not
lie in $\overline{Y}$, while every other $T$-orbit does contain a
point of $\overline{Y}$.  
\end{example}

A fundamental problem in algebraic geometry, especially in the
consideration of moduli spaces, is to find a good compactification of
a variety.  Here the notion of ``good'' will depend on the specific
problem, but we often require the compactification to be smooth,
the ``boundary'' (new points added in the compactification) to be a
divisor (codimension-one), and the irreducible components of this
divisor to intersect nicely.

\begin{example} \label{e:DeconciniProcesi}
Let $\mathcal A = \{ H_1,\dots,H_s \}$ be a hyperplane arrangement in
$\mathbb P^n$, where $H_i = \{ x \in \mathbb P^n : a_i \cdot x =0\}$,
with $a_i \in \mathbb C^{n+1}$.  Let $Y = \mathbb P^n \setminus
\mathcal A$.  This can be embedded into $(\mathbb C^*)^{s-1}$ by
sending $y \in Y$ to $[a_1 \cdot y: \dots: a_s \cdot y]$.  Let $A$ be
the $(n+1) \times s$ matrix with columns the vectors $a_i$, which we
assume has rank $n+1$, and let $B$ be a $(s-n-1) \times (n+1)$ matrix
of rank $s-n-1$ with $AB^{T}=0$.  Then $Y = V(\sum_{i=1}^n b_{ij}x_j :
1 \leq i \leq s-n-1)$.  One choice of compactification of $Y$ is the
original $\mathbb P^n$; another is the DeConcini-Procesi wonderful
compactification~\cite{DeConciniProcesi}.
\end{example}

A natural way to compactify a variety is to take its closure in a
larger variety.  The idea we follow, due to
Tevelev~\cite{TevelevCompactifications}, is to take as the larger
variety a toric variety whose fan has support the tropical variety.
Note that if $\dim(Y)>1$ there is not a unique such fan, and there may
be no preferred choice.

\begin{definition}
Fix $Y \subset T$, and choose a fan $\Sigma$ with support
$\trop(Y)$.  The closure $\overline{Y} = \cl(Y \subset
X_{\Sigma})$ is a {\em tropical compactification} of $Y$.
\end{definition}

Tropical compactifications have nice properties:
\begin{enumerate}
\item $\overline{Y}$ is proper (compact),
\item $\overline{Y}$ intersects a codimension-$k$ $T$-orbit of
  $X_{\Sigma}$ in codimension $k$.
\end{enumerate}

The first condition is somewhat surprising, since the toric variety
$X_{\Sigma}$ is not itself compact.  For example, if $Y = V(x+y+1)
\subseteq (\mathbb C^*)^2$, then $\trop(Y)$ is the standard tropical
line.  There is only one polyhedral fan $\Sigma$ with this support,
and the corresponding toric variety is $\mathbb P^2$ with three points
removed, which is not compact.  However in this case the closure of
$Y$ in $\mathbb P^2$, which is clearly compact, does not pass through
these three points, so the closure of $Y$ in $X_{\Sigma}$ equals the
closure of $Y$ in $\mathbb P^2$, so is also compact.

The second condition means that the pull-backs to $\overline{Y}$ of
the codimension-$k$ torus-invariant strata on $X_{\Sigma}$ give a
finite number of distinguished classes in the Chow group $A^k(\overline{Y})$.  For
example when $k=1$ this gives a distinguished collection of effective divisors on
$\overline{Y}$.

If the fan $\Sigma$ is chosen to be sufficiently refined then we get
further control over the compactification $\overline{Y}$.  One way to
guarantee ``sufficiently refined'' is to choose a fan $\Sigma$ so that
$\inn_w(I(Y))$ is constant for all $w \in \relint{\sigma}$ for all
$\sigma \in \Sigma$.  By further refining $\Sigma$ we can also assume
that the toric variety $X_{\Sigma}$ is smooth.  With a sufficiently
refined fan we have:

\begin{enumerate}
\item The multiplicity $w_{\sigma}$ on a maximal cone $\sigma$ of
  $\Sigma$ equals the intersection number $[\overline{Y}] \cdot
  [V(\sigma)]$, where $V(\sigma)$ is the closure of the $T$-orbit
  corresponding to $\sigma$.  
\item For any $\sigma \in \Sigma$, the intersection of $\overline{Y}$
  with the $T$-orbit $\mathcal O_{\sigma}$ is the quotient $V(\inn_w(I(Y)))/(\mathbb C^*)^{\dim
    \sigma}$ for any $w \in \mathrm{relint}(\sigma)$.
\end{enumerate}

See also \cite{Hacking}, \cite{SpeyerThesis},  and
\cite{SturmfelsTevelev} for more details on this.  The connection to
intersection theory has been further developed by Allermann and
Rau~\cite{AllermannRau}, which is the tropical analogue of earlier
toric work by Fulton and Sturmfels~\cite{FultonSturmfels}.

\begin{example}
When $Y = \mathbb P^n \setminus \mathcal A$, as in Example~\ref{e:DeconciniProcesi}, there is a
coarsest fan structure $\Sigma$ on $\trop(Y) \subseteq \mathbb
R^{s-1}$, since every other fan with support $\trop(Y)$ has cones that
subdivide those of $\Sigma$. The tropical compactification
$\overline{Y}$ coming from taking the closure in $X_{\Sigma}$ is the
DeConcini/Procesi wonderful compactification for most choices of
$\mathcal A$; see~\cite{FeichtnerSturmfels} and \cite{TevelevCompactifications}.
\end{example}

 A motivating example of this is given by the moduli space
 $\overline{M}_{0,n}$.  The moduli space $M_{0,n}$ parameterizes
 isomorphism classes of smooth genus zero curves with $n$ distinct
 marked points.  It thus parameterizes ways to arrange $n$ distinct
 point on $\mathbb P^1$ up to $\Aut(\mathbb P^1)$.  For example,
 $M_{0,3}$ is a point, as there is an automorphism of $\mathbb P^1$
 that takes any three distinct points to $0, 1, \infty$.  When $n=4$,
 $M_{0,4}=\mathbb P^1 \setminus \{ 0, 1, \infty \}$.  In general,
\begin{align*}
M_{0,n} & =(\mathbb P^1 \setminus \{0,1, \infty \})^{n-3} \setminus
\text{ diagonals}\\ & = (\mathbb C^* \setminus \{1 \})^{n-3} \setminus
\text{ diagonals}\\ & = \mathbb P^{n-3} \setminus \{ x_0=0, x_i=0,
x_i=x_0, x_i=x_j : 1 \leq i<j \leq n \}. \\
\end{align*}
We thus have $M_{0,n}$ as the complement of ${n -1 \choose 2} = {n
  \choose 2} -n +1$ hyperplanes.  This means the moduli space
$M_{0,n}$ can be be embedded into $(\mathbb C^*)^{{n \choose 2}-n}$ as
a closed subvariety as in Example~\ref{e:DeconciniProcesi}.  The
tropical variety $\trop(M_{0,n})$ is the support of a  $(n-3)$-dimensional fan
$\Delta$ in $\mathbb R^{{n \choose 2}-n}$.  The toric variety
$X_{\Delta}$ is smooth, but not complete (projective).

The fan $\Delta$ is the space of phylogenetic trees, which comes from
mathematical biology.  See, for example, \cite{BilleraHolmesVogtmann}
for more on this space, and \cite{SpeyerSturmfelsTropicalGrassmannian}
for some of the connection, which uses Kapranov's
description~\cite{KapranovChow} of $\MOn$ as the Chow quotient of a
Grassmannian.  Maximal cones of $\Delta$ are labelled by trivalent trees with $n$
labelled leaves.  A point in the cone records the length of the
internal edges in the tree.

\begin{figure}
\center{\resizebox{!}{9cm}{\input{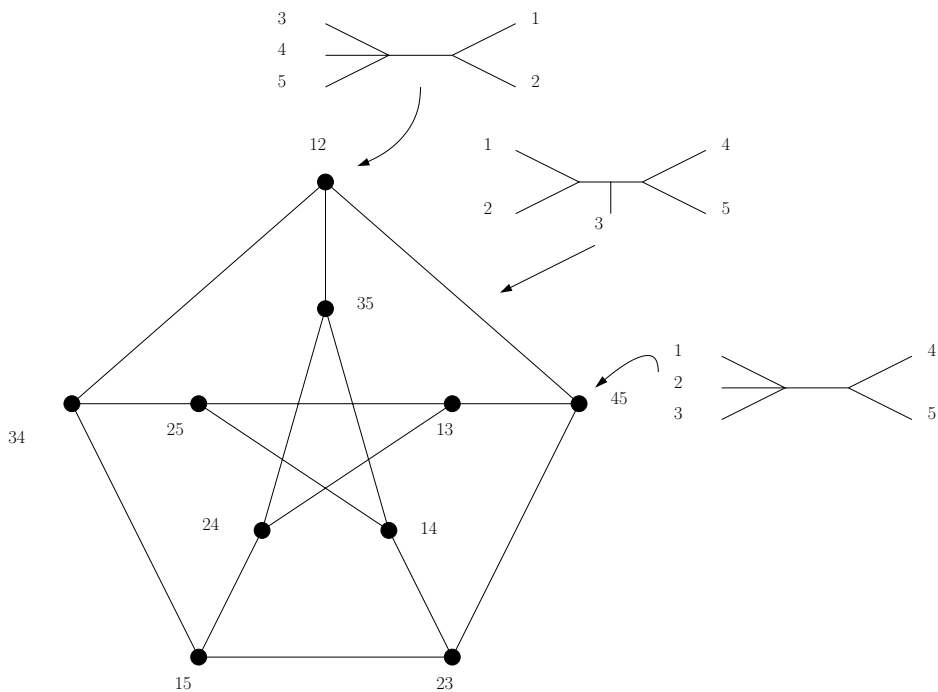}}}
\caption{The tropical variety of $M_{0,5}$ \label{f:M05}}
\end{figure}

A picture of $\Delta$ when $n=5$ is shown in Figure~\ref{f:M05}.  This
is a two-dimensional fan in $\mathbb R^5$, so its intersection with
the four-dimensional sphere in $\mathbb R^5$ is a graph, which is
drawn in Figure~\ref{f:M05}.

The closure of $M_{0,n}$ in $X_{\Delta}$ is the Deligne-Mumford moduli
space $\overline{M}_{0,n}$ of {\em stable} genus zero curves with $n$
marked points.  This parameterizes isomorphism classes of trees of
$\mathbb P^1$'s with $n$ marked points and at least three special
points (nodes or marked points) on each component.
See~\cite{GibneyMaclagan} or \cite{TevelevCompactifications} for
details.

\begin{figure}
\center{\resizebox{!}{5cm}{\input{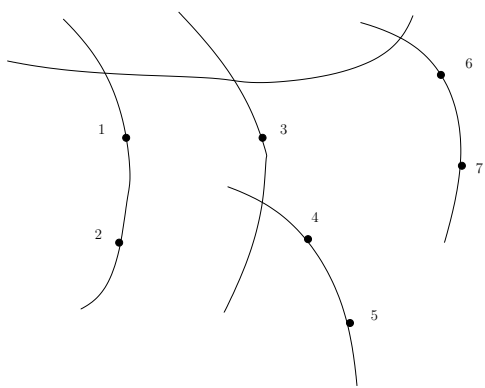}}}
\caption{\label{f:stablecurve} A stable curve with $7$ marked points}
\end{figure}

 For $\sigma \in \Delta$, the intersection of
$\overline{M}_{0,n}$ with the torus orbit corresponding to $\sigma$ is
the stratum of all curves with dual graph the corresponding tree.  In
particular, the intersection of $\overline{M}_{0,n}$ with a
torus-invariant divisor on $X_{\Delta}$ is a boundary divisor.

The moduli space $\overline{M}_{0,n}$ and the toric variety
$X_{\Delta}$ are closely related.  Their Picard groups are isomorphic,
and the inclusion $i: \overline{M}_{0,n} \rightarrow X_{\Delta}$
introduces an isomorphism $i^* : A^*(X_{\Delta}) \rightarrow
A^*(\overline{M}_{0,n})$.  This generalizes to general wonderful
compactifications; see~\cite{FeichtnerYuzvinsky}.
\bibliographystyle{amsplain}

\def\cprime{$'$}
\providecommand{\bysame}{\leavevmode\hbox to3em{\hrulefill}\thinspace}
\providecommand{\MR}{\relax\ifhmode\unskip\space\fi MR }
\providecommand{\MRhref}[2]{%
  \href{http://www.ams.org/mathscinet-getitem?mr=#1}{#2}
}
\providecommand{\href}[2]{#2}

\end{document}